\documentclass{amsart}
\usepackage{amscd,amssymb,graphicx}
\author[Bodin]{Derek Bodin}
\address{
Derek Bodin\\
University of Minnesota\\
Minneapolis, MN 55455 }
\email{bodin@cs.umn.edu}
\author[DeCleene]{Chris DeCleene}
\address{
Chris DeCleene\\
University of Wisconsin-Eau Claire\\
Eau Claire, WI 54702-4004}
\email{cdecleene@gmail.com}
\author[Hager]{William Hager}
\address{
William Hager\\
University of Iowa\\
Iowa City, IA 52245-4027}
\email{william-hager@uiowa.edu}
\author[Otto]{Carolyn Otto}
\address{
Carolyn Otto\\
Rice University\\
Houston, TX 77005-1827} \email{cotto@rice.edu}
\author[Penkava]{Michael Penkava}
\address{
Michael Penkava\\
University of Wisconsin-Eau Claire\\
Eau Claire, WI 54702-4004} \email{penkavmr@uwec.edu}
\author[Phillipson]{Mitch Phillipson}
\address{
Mitch Phillipson\\
University of Wisconsin-Eau Claire\\
Eau Claire, WI 54702-4004}
\email{phillima@uwec.edu}
\author[Steinbach]{Ryan Steinbach}
\address{
Ryan Steinbach\\
University of Wisconsin-Madison\\
Madison, WI 53706-1796}
\email{rsteinbach@wisc.edu}
\author[Weber]{Eric Weber}
\address{
Eric Weber\\
University of Wisconsin-Eau Claire\\
Eau Claire, WI 54702-4004}
\email{webered@uwec.edu}
\thanks{Research of these authors was partially supported by grants
from the National Science Foundation and the University of
Wisconsin-Eau Claire.}

\newtheorem{thm}{Theorem}[section]
\newtheorem{con}[thm]{Conjecture}

\theoremstyle{definition}

\newtheorem{rem}{Remark}

\def \ph{\varphi}

\def\Ch{\operatorname{Ch}}

\def \ra{\rightarrow}

\def \hom{\mbox{\rm Hom}}
\def \ie{\hbox{\it i.e.}}

\def \tns{\otimes}

\def \mplus{+\cdots+}
\def \mcom{,\cdots,}
\def \k{\mbox{$\mathbb K$}}

\def \C{\mbox{$\mathbb C$}}
\def \Z{\mbox{$\mathbb Z$}}

\def\br#1#2{\lbrack#1,#2\rbrack}

\def\zt{\mbox{$\Z_2$}}

\def\inv{^{-1}}

\def\im{\operatorname{Im}}

\def\m{\mbox{$\mathfrak m$}}

\def\coder{\operatorname{Coder}}
\def\ainf{\mbox{$A_\infty$}}

\def\and{\mbox{ \rm and }}
\def\T{\mathcal T}
\def\TV{\T(V)}
\def\TW{\mbox{$\T(W)$}}

\def\s#1{(-1)^{#1}}
\DeclareMathOperator*{\invlim}{\underleftarrow{\rm lim}}

\def\pha#1#2{\ph^{#1}_{#2}}

\def\psa#1#2{\psi^{#1}_{#2}}

\def\thmref#1{Theorem (\ref{#1})}
\def\inv{^{-1}}

\def\dinf{\mbox{$d^\text{inf}$}}
\def\P{\mathbb P}

\newcommand{\mapsiso}{\tilde{\ra}}

\def\dec{\operatorname{Ch}}
\def\Cn#1#2{C^{#1}_{#2}}

\input{xy}
\xyoption{all}

\begin{document}
\setlength{\multlinegap}{0pt}
\title[The moduli space of $1|1$-dimensional complex associative algebras]
{Moduli space of $1|1$-dimensional complex associative algebras}%

\address{}%
\email{}%

\thanks{}%

\date{\today}
\begin{abstract}
In this paper, we study the moduli space of $1|1$-dimensional complex
associative algebras. We give a complete calculation of the cohomology of every
element in the moduli space, as well as compute their versal deformations.

\end{abstract}
\maketitle

\nocite{ps1,ps2,mp1,pv1}
\section{Introduction}

Super Lie algebras, or \zt-graded Lie algebras, have been studied for a long
time, and have many applications in mathematics and physics. The notion of a
\zt-graded associative algebra is not as well known (however, see \cite{knus}), but these algebras are
examples of \zt-graded \ainf\ algebras, and thus they arise naturally in the
study of \ainf\ algebras.  Although we will not consider extensions of
\zt-graded associative algebras to more general \ainf\ structures in this
paper, the results here are the first step in the construction of such
extensions. We plan to discuss such extensions in a later paper, restricting
ourselves here to giving a complete description of the moduli space of
$1|1$-dimensional associative algebras.

In the case of Lie algebras, the \zt-graded Jacobi identity picks up some signs
that depend on the parity of the elements being bracketed, but the
associativity relation for \zt-graded associative algebras does not pick up any
signs, so it may seem at first glance that there are no new features which
arise in the study of \zt-graded associative algebras.

The moduli space of equivalence classes of \zt-graded associative algebras on a
vector space of dimension $m|n$ differs from the moduli space of associative
structures on the same space ignoring the grading in two important ways. First,
a \zt-graded algebra structure is required to be an even map, which means that
not all associative algebra structures are allowed in the \zt-graded case.
Secondly, the moduli space is given by equivalence classes of algebra
structures under an action by the group of linear automorphisms of the vector
space. 

For the \zt-graded case, we only allow even automorphisms, which means
that the equivalence classes are potentially smaller in the \zt-graded case.
Since there are fewer allowable \zt-graded algebra structures, but also fewer
equivalences between them, it is not obvious whether the moduli space of
\zt-graded associative algebras on a \zt-graded vector space is larger or
smaller than the moduli space of all associative algebra structures on the
vector space.

There is a map between the moduli space of \zt-graded
algebra structures on a \zt-graded  vector space and the moduli space of all
algebra structures on the underlying space. This map in general is neither
injective nor surjective. In fact, there are exactly 6 isomorphism classes of
$1|1$-dimensional associative algebras, and exactly 6 isomorphism classes of
ordinary associative algebras on a 2-dimensional vector space, but the map
between the \zt-graded algebras to the ordinary ones has exactly 5 algebras in
its image, so that two of the \zt-graded associative algebras map to the same
image. Thus, even in the simplest case, the map between the moduli spaces
is neither injective nor surjective.

One method of constructing the moduli space of ordinary associative algebras in
dimension 2 is to consider extensions of a 1-dimensional associative algebra by
a 1-dimensional associative algebra. This is possible because there are no
simple 2-dimensional associative algebras, by a theorem of Wedderburn, so all
such algebras have an ideal, and therefore arise as extensions.  As we will see
in this paper, there is a simple $1|1$-dimensional associative algebra, so the
theorem of Wedderburn, in its classical form, does not apply in the \zt-graded
case.  Therefore, we will use a different method of determining the equivalence
classes in this paper. However, there is a natural generalization of Wedderburn's
theorem to \zt-graded algebras, and the simple $1|1$-dimensional algebra plays an
important role in this generalization, because it is a \zt-graded division algebra
over \C.  We will not discuss this issue further in this paper, but refer the reader
to \cite{Karrer,bppw1,bsz}.

In this paper, we will give a complete description of the moduli space of
$1|1$-dimensional algebras, including a computation of a miniversal deformation
of each of these algebras.  From the miniversal deformations, a decomposition
of the moduli space into strata is obtained, with the only connections between
strata given by jump deformations. In the $1|1$-dimensional case, the
description is simple, because each of the strata consists of a single point, so
the only interesting information is given by the jump deformations.

The versal deformation of an associative algebra depends only on the second and
third Hochschild cohomology groups.  However, we give a complete calculation of
the cohomology for each of the algebras. What makes the study of associative
algebras of low dimension much more complicated than the corresponding study of
low dimensional Lie algebras is that while for a Lie algebra, the $n$-th cohomology
group $H^n$ vanishes
for $n$ larger than the dimension of the vector space, in general, for an
associative algebra $H^n$ does not vanish.  Thus we had to develop arguments on
a case by case basis for each of the six distinct algebras.  In particular, one
of these algebras has an unusual pattern for the cohomology, which made its
computation rather nontrivial.

The main result of this paper is the complete description of the Hochschild 
cohomology for
all $1|1$-dimensional associative algebras. It turns out that the calculation
of cohomology even for low dimensional associative algebras is a nontrivial
problem. To construct extensions of associative algebras to \ainf\ algebras, it
is necessary to have a complete description of the cohomology in all degrees,
not just $H^2$ and $H^3$, which are needed for the deformation theory of these
algebras as associative algebras. What we compute in this paper is the first
step in constructing $1|1$-dimensional \ainf\ algebras. These results may be of
interest on their own, especially as an indication of the difficulty which
occurs in computing the deformation theory of associative algebras, even in low
dimension.

\section{Preliminaries}
Suppose that $V$ is a vector space, defined over a field $\k$
whose characteristic is not 2 or 3, equipped with an associative multiplication
structure $m:V\tns V\ra V$. The associativity relation can be given in the form
\begin{equation*}
m\circ(m\tns 1)=m\circ(1\tns m).
\end{equation*}

When the space $V$ is \zt-graded, there is no difference in the relation of
associativity, but only even maps $m$ are allowed, so the set of associative
algebra structures depends on the \zt-grading in this way.

The notion of \emph{equivalence} of associative algebra structures is given as
follows. If $g$ is a linear automorphism of $V$, then define
\begin{equation*}
g^*(m)=g\inv\circ m\circ (g\tns g).
\end{equation*}
Two algebra structures $m$ and $m'$ are equivalent if there is an automorphism
$g$ such that $m'=g^*(m)$. The set of equivalence classes of algebra structures
on $V$ is called the \emph{moduli space} of associative algebras on $V$.

When $V$ is \zt-graded, we require that $g$ be an even map. Thus the set of
equivalence classes of \zt-graded associative algebra structures will be
different than the set of equivalence classes of associative algebra structures
on the same space, ignoring the grading. Because the set of equivalences is
more restricted in the \zt-graded case, two algebra structures which are
equivalent as ungraded algebra structures may not be equivalent as \zt-graded
algebra structures. There is a map between the moduli space of \zt-graded
algebra structures on $V$ and the space of all algebra structures on $V$. In
general, this map will be neither injective nor surjective.

\emph{Hochschild cohomology} was introduced in \cite{hoch}, and used to
classify infinitesimal deformations of associative algebras. Suppose that
\begin{equation*}
m_t=m+t\ph,
\end{equation*}
is an infinitesimal deformation of $m$.  By this we mean that the structure
$m_t$ is associative up to first order. From an algebraic point of view, this
means that we assume that $t^2=0$, and then check whether associativity holds.
It is not difficult to show that is equivalent to the following.
\begin{equation*}
a\ph(b,c)-\ph(ab,c)+\ph(a,bc)-\ph(a,b)c=0,
\end{equation*}
where, for simplicity, we denote $m(a,b)=ab$. Moreover, if we let
\begin{equation*}
g_t=I+t\lambda
\end{equation*} be an infinitesimal automorphism of $V$, where
$\lambda\in\hom(V,V)$, then it is easily checked that
\begin{equation*}
g_t^*(m)(a,b)=ab+t(a\lambda(b)-\lambda(ab)+\lambda(a)b).
\end{equation*}
 This naturally leads
to a definition of the Hochschild coboundary operator $D$ on $\hom(\TV,V)$ by
\begin{align*}
D(\ph)(a_0\mcom a_n)=&a_0\ph(a_1\mcom a_n)+\s{n+1}\ph(a_0\mcom a_{n-1})a_n\\&+\sum_{i=0}^{n-1}\s{i+1}\ph(a_0\mcom a_{i-1},a_ia_{i+1},a_{i+2}\mcom a_n)
.
\end{align*}
If we set $C^n(V)=\hom(V^n,V)$, then $D:C^n(V)\ra C^{n+1}(V)$. One obtains the
following classification theorem for infinitesimal deformations.
\begin{thm} The equivalence classes of infinitesimal deformations $m_t$ of an associative algebra structure $m$ under the action of the
group of infinitesimal automorphisms on the set of infinitesimal deformations
are classified by the Hochschild cohomology group
\begin{equation*}
H^2(m)=\ker(D:C^2(V)\ra C^3(V))/\im(D:C^1(V)\ra C^2(V)).
\end{equation*}
\end{thm}
When $V$ is \zt-graded, the  only modifications that are necessary are that
$\ph$ and $\lambda$ are required to be even maps, so we obtain that the
classification is given by $H^2_e(V)$, the even part of the Hochschild
cohomology.

We wish to transform this classical viewpoint into the more modern viewpoint of
associative algebras as being given by codifferentials on a certain coalgebra.
To do this, we first introduce the \emph{parity reversion} $\Pi V$ of a
\zt-graded vector space $V$. If $V=V_e\oplus V_o$ is the decomposition of $V$
into its even and odd parts, then $W=\Pi V$ is the \zt-graded vector space
given by $W_e=V_o$ and $W_o=V_e$. In other words, $W$ is just the space $V$
with the parity of elements reversed.

Denote the tensor (co)-algebra of $W$ by $\TW=\bigoplus_{k=0}^\infty W^k$,
where $W^k$ is the $k$-th tensor power of $W$ and $W^0=\k$. For brevity, the
element in $W^k$ given by the tensor product of the elements $w_i$ in $W$ will
be denoted by $w_1\cdots w_k$. The coalgebra structure on $\TW$ is  given by
\begin{equation*}
\Delta(w_1\cdots w_n)=\sum_{i=0}^n w_1\cdots w_i\tns w_{i+1}\cdots w_n.
\end{equation*}
Define $d:W^2\ra W$ by $d=\pi\circ m\circ (\pi\inv\tns\pi\inv)$, where
$\pi:V\ra W$ is the identity map, which is odd, because it reverses the parity
of elements. Note that $d$ is an odd map. The space $C(W)=\hom(\TW,W)$ is
naturally identifiable with the space of coderivations of $\TW$.  In fact, if
$\ph\in C^k(W)=\hom(W^k,W)$, then $\ph$ is extended to a coderivation of $\TW$
by
\begin{equation*}
\ph(w_1\cdots w_n)=
\sum_{i=0}^{n-k}\s{(w_1\mplus w_i)\ph}w_1\cdots
 w_i\ph(w_{i+1}\cdots w_{i+k})w_{i+k+1}\cdots w_n.
\end{equation*}

The space of coderivations of $\TW$ is equipped with a \zt-graded Lie algebra
structure given by
\begin{equation*}
[\ph,\psi]=\ph\circ\psi-\s{\ph\psi}\psi\circ\ph.
\end{equation*}
The reason that it is more convenient to work with the structure $d$ on $W$
rather than $m$ on $V$ is that the condition of associativity for $m$
translates into the codifferential property $[d,d]=0$.  Moreover, the
Hochschild coboundary operation translates into the coboundary operator $D$ on
$C(W)$, given by
\begin{equation*}
D(\ph)=[d,\ph].
\end{equation*}
This point of view on Hochschild cohomology first appeared in \cite{sta4}.  The
fact that the space of Hochschild cochains is equipped with a graded Lie
algebra structure was noticed much earlier \cite{gers,gers1,gers2,gers3,gers4}.

For notational purposes, we introduce a basis of $C^n(W)$ as follows.  Suppose
that $W=\langle w_1\mcom w_m\rangle$. Then if $I=(i_1\mcom i_n)$ is a
\emph{multi-index}, where $1\le i_k\le m$, denote $w_I=w_{i_1}\cdots w_{i_n}$.
Define $\ph^{I}_i\in C^n(W)$ by
\begin{equation*}
\ph^I_i(w_J)=\delta^I_Jw_i,
\end{equation*}
where $\delta^I_J$ is the Kronecker delta symbol. In order to emphasize the
parity of the element, we will denote $\ph^I_i$ by $\psi^I_i$ when it is an odd
coderivation.

For a multi-index $I=(i_1\mcom i_k)$, denote its \emph{length}  by $\ell(I)=k$.  If
$K$ and $L$ are multi-indices, then denote $KL=(k_1\mcom k_{\ell(K)},l_l\mcom
l_{\ell(L)})$.  Then
\begin{align*}
(\ph^I_i\circ\ph^J_j)(w_K)&=
\sum_{K_1K_2K_3=K}\s{w_{K_1}\ph^J_j} \ph^I_i(w_{K_1},\ph^J_j(w_{K_2}), w_{K_3})
\\&=
\sum_{K_1K_2K_3=K}\s{w_{K_1}\ph^J_j}\delta^I_{K_1jK_3}\delta^J_{K_2}w_i,
\end{align*}
from which it follows that
\begin{equation}\label{braform}
\ph^I_i\circ\ph^J_j=\sum_{k=1}^{\ell(I)}\s{(w_{i_1}\mplus w_{i_{k-1}})\ph^J_j}
\delta^k_j
\ph^{(I,J,k)}_i,
\end{equation}
where $(I,J,k)$ is given by inserting $J$ into $I$ in place of the $k$-th
element of $I$; \ie, $(I,J,k)=(i_1\mcom i_{k-1},j_1\mcom j_{\ell(J)},i_{k+1}\mcom
i_{\ell(I)})$.

Let us recast the notion of an infinitesimal deformation in terms of the
language of coderivations.  We say that
\begin{equation*}
d_t=d+t\psi
\end{equation*}
is a deformation of the codifferential $d$ precisely when $[d_t,d_t]=0 \mod t^2$.
This condition immediately reduces to the cocycle condition $D(\psi)=0$.  Note
that we require $d_t$ to be odd, so that $\psi$ must be an odd coderivation.
One can introduce a more general idea of parameters, allowing both even and odd
parameters, in which case even coderivations play an equal role, but we will
not adopt that point of view in this paper.

For associative algebras, we require that $d$ and $\psi$ lie in $\hom(W^2,W)$.
This notion naturally generalizes to considering $d$ simply to be an arbitrary
odd codifferential, in which case we would obtain an \ainf\ algebra, a natural
generalization of an associative algebra.

More generally, we need the notion of a versal deformation, in order to understand how the moduli space is glued together. To explain versal deformations we introduce the notion of  a deformation with a local base.

A local base $A$ is a \zt-graded commutative, unital $\k$-algebra with an augmentation $\epsilon:A\ra\k$, whose kernel $\m$ is the unique maximal ideal in $A$, so that $A$ is a local ring. It follows that $A$ has a unique decomposition $A=\k\oplus\m$ and $\epsilon$ is just the projection onto the first factor. Let $W_A=W\tns A$ equipped with the usual structure of a right $A$-module. Let $T_A(W_A)$ be the tensor algebra of $W_A$ over $A$, that is $T_A(W_A)=\bigoplus_{k=0}^\infty T^k_A(W_A)$ where $T^0_A(W_A)=A$ and $T^{k+1}_A(W_A)=T^k(W_A)_A\tns_A W_A$. It is a standard fact that $T^k_A(W_A)=T^k(W)\tns A$ in a natural manner, and thus $T_A(W_A)=T(W)\tns A$.

Any $A$-linear map $f:T_A(W)\ra T_A(W)$ is induced by its restriction to $T(W)\tns \k=T(W)$ so we can view an $A$-linear coderivation $\delta_A$ on $T_A(W_A)$ as a map $\delta_A:T(W)\ra T(W)\tns A$. A morphism $f:A\ra B$ induces a map  $$f_*:\coder_A(T_A(W_A))\ra \coder_B(T_B(W_B))$$ given by $f_*(\delta_A)=(1\tns f)\delta_A$, moreover if $\delta_A$ is a codifferential then so is $f_*(A)$. A codifferential $d_A$ on $T_A(W_A)$ is said to be a deformation of the codifferential $d$ on $T(W)$ if $\epsilon_*(d_A)=d$.

If $d_A$ is a deformation of $d$ with base $A$ then we can express
\begin{equation*}
 d_A=d+\ph
\end{equation*}
where $\ph:T(W)\ra T(W)\tns\m$. The condition for $d_A$ to be a codifferential is the Maurer-Cartan equation,
\begin{equation*}
D(\ph)+\frac12[\ph,\ph]=0
\end{equation*}
If $\m^2=0$ we say that $A$ is an infinitesimal algebra and a deformation with base $A$ is called infinitesimal.

A typical example of an infinitesimal base is $\k[t]/(t^2)$, moreover, the classical notion of an infinitesimal deformation
$$d_t=d+t\ph$$
is precisely an infinitesimal deformation with base $\k[t]/(t^2)$.

A local algebra $A$ is complete if
\begin{equation*}
 A=\invlim_kA/\m^k
\end{equation*}
A complete, local augmented $\k$-algebra will be called formal and a deformation with a formal base is called a formal deformation. An infinitesimal base is automatically formal, so every infinitesimal deformation is a formal deformation.

An example of a formal base is $A=\k[[t]]$ and a deformation of $d$ with base $A$ can be expressed in the form
$$d_t=d+t\psi_1+t^2\psi_2+\dots$$
This is the classical notion of a formal deformation. It is easy to see that the condition for $d_t$ to be a formal deformation reduces to
\begin{align*}
 D(\psi_{n+1})=-\frac12\sum_{k=1}^{n}[\psi_k,\psi_{n+1-k}]
\end{align*}

An automorphism of $W_A$ over $A$ is an $A$-linear isomorphism $g_A:W_A\ra W_A$ making the diagram below commute.
\begin{figure}[h!]
 $$\xymatrix{
 W_A \ar[r]^{g_A} \ar[d]^{\epsilon_*} & W_A \ar[d]^{\epsilon_*} \\
 W \ar[r]^I & W}$$
\end{figure}
The map $g_A$ is induced by its restriction to $T(W)\tns\k$ so we can view $g_A$ as a map
$$g_A:T(W)\ra T(W)\tns A$$
so we ca express $g_A$ in the form
$$g_A=I+\lambda$$
where $\lambda:T(W)\ra T(W)\tns\m$. If $A$ is infinitesimal then $g_A^{-1}=I-\lambda$.

Two deformations $d_A$ and $d_A'$ are said to be equivalent over $A$ if there is an automorphism $g_A$ of $W_A$ over $A$ such that $g_A^*(d_A)=d_A'$. In this case we write $d'_A\sim d_A$.

An infinitesimal deformation $d_A$ with base $A$ is called universal if whenever $d_B$ is an infinitesimal deformation with base $B$, there is a unique morphism $f:A\ra B$ such that $f_*(d_A)\sim d_B$.

\begin{thm}
 If $\dim H^2_{odd}(d)<\infty$ then there is a universal infinitesimal deformation $\dinf$ of $d$. Given by
 $$\dinf=d+\delta^it_i$$
 where $H^2_{odd}(d)=\langle\bar{\delta^i}\rangle$ and $A=\k[t_i]/(t_it_j)$ is the base of deformation.
\end{thm}

A formal deformation $d_A$ with base $A$ is called versal if given any formal deformation of $d_B$ with base $B$ there is a morphism $f:A\ra B$ such that $f_*(d_A)\sim d_B$. Notice that the difference between the versal and the universal property of infinitesimal deformations is that $f$ need not be unique. A versal deformation is called \emph{miniversal} if $f$ is unique whenever $B$ is infinitesimal. The basic result about versal deformation is:
\begin{thm}
 If $\dim H^2_{odd}(d)<\infty$ then a miniversal deformation of $d$ exists.
\end{thm}

In this paper we will only need the following result to compute the versal deformations.
\begin{thm}
 Suppose $H^2_{odd}(d)=\langle\bar{\delta^i}\rangle$ and $[\delta^i,\delta^j]=0$ for all $i,j$ then the infinitesimal deformation
 $$\dinf=d+\delta^it_i$$
 is miniversal, with base $A=\k[[t_i]].$
\end{thm}

The construction of the moduli space as a geometric object is based on the idea that codifferentials
which can be obtained by deformations with small parameters are ``close'' to each other. From the small deformations,
we can construct 1-parameter families or even multi-parameter families, which are defined for small values of the
parameters, except possibly when the parameters vanish.

If $d_t$ is a one parameter family of deformations, then two things can occur. First, it may happen that
$d_t$ is equivalent to a certain codifferential $d'$ for every small value of $t$ except zero.
Then we say that $d_t$ is a jump deformation from $d$ to $d'$. It will
never occur that $d'$ is equivalent to $d$, so  there are no jump deformations from a codifferential to itself.
Otherwise, the codifferentials $d_t$ will all be nonequivalent if $t$ is small enough. In this case, we
say that $d_t$ is a smooth deformation.

In \cite{fp10}, it was proved for Lie algebras that given three codifferentials $d$, $d'$ and $d''$,
if there are jump deformations from $d$ to $d'$ and from $d'$ to $d''$, then  there is a jump deformation from
$d$ to $d''$. The proof of the corresponding statement for associative algebras is essentially the same.

Similarly, if there is a jump deformation from $d$ to $d'$, and a family of smooth deformations
$d'_t$, then there is a family $d_t$ of smooth deformations of $d$, such that every deformation in the image of
$d'_t$ lies in the image of $d_t$, for sufficiently small values of $t$. In this case, we say that the
smooth deformation of $d$ factors through the jump deformation to $d'$.

In the examples of complex moduli spaces of Lie and associative algebras which we have studied,
it turns out that there is a natural
stratification of the moduli space of $n$-dimensional algebras by orbifolds, where the codifferentials
on a given strata are connected by smooth deformations, which don't factor through jump deformations. These smooth deformations determine the local neighborhood structure.

The strata are connected by jump deformations, in the sense that
any smooth deformation from a codifferential on one strata
to another strata factors through a jump deformation.  Moreover, all of the strata are given by projective
orbifolds. In fact, in all the complex examples we have studied, the orbifolds either are single points,
or $\C\P^n$ quotiented out by either $\Sigma_{n+1}$ or a subgroup, acting on
 $\C\P^n$ by permuting the coordinates.

We don't have a concrete proof at this time, but we conjecture that this pattern holds in general. In other words, we believe the following conjecture.
\begin{con}
The moduli space of Lie or associative algebras of a fixed finite dimension $n$ are stratified by projective orbifolds,
with jump deformations and smooth deformations factoring through jump deformations providing the only
deformations between the strata.
\end{con}

\section{Associative algebra structures on a $1|1$ vector space}
Suppose that $W=\langle e,f \rangle$,  where $e$ is an even element, and $f$ is
odd.  Then $C^n=\langle \ph^I_i,\ell(I)=n\rangle$ has dimension $\dim C^n=2^{n+1}$.
For later convenience, we decompose $C^n$ as follows. Let
\begin{align*}
E^n=\langle \ph^I_e, \ell(I)=n\rangle\\
F^n=\langle \ph^I_f,\ell(I)=n\rangle.
\end{align*}
Then $C^n=E^n\oplus F^n$. Moreover $\dim E^n=\dim F^n=2^n$.

Now, an odd element in $C^2$ is of the form $d=\psa{ff}fx+\psa{fe}ey+\psa{ef}ez
+\psa{ee}fw$. One computes that
\begin{align*}
\tfrac12[d,d]=&\pha{eff}ey(x+y)+\pha{eee}ew(y+z)+\pha{efe}fw(y+z)\\&+\pha{eef}fw(x+y)
+\pha{ffe}ez(x-z)-\pha{fee}fw(x-z).
\end{align*}
Setting $[d,d]=0$, we obtain 6 distinct, nonequivalent codifferentials.
\begin{align*}
d_1&=-\psi_f^{ff}+\psi_e^{ef}-\psi_e^{fe}+\psi_f^{ee}\\
d_2&=\psi_f^{ff}\\
d_3&=-\psi_f^{ff}+\psi_e^{ef}\\
d_4&=\psi_f^{ff}+\psi_e^{fe}\\
d_5&=-\psi_f^{ff}+\psi_e^{ef}-\psi_e^{fe}\\
d_6&=\psi_f^{ee}.
\end{align*}
 Note that if we
define $D(\ph)=\br {d^*}\ph$, where $d^*$ is one of the above codifferentials,
then $D^2=0$, so the \emph{coboundary operator} $D$ determines a differential
on $C(W)$. Since $d^*\in C^2$, $D(C^k)\subseteq C^{k+1}$, and we can define the
$k$-th cohomology $H^k(d^*)$ of $d^*$ by
\begin{equation*}
H^k(d^*)=\ker(d^*:C^k\ra
C^{k+1})/\im(d^*:C^{k-1}\ra C^k).
\end{equation*}

The cohomology of these codifferentials is given in Table \ref{coho table}
below. These codifferentials can be distinguished in terms of their cohomology,
with the exception of $d_2$ and $d_3$, which are opposite algebras.

\begin{table}[h!]
\begin{center}
\begin{tabular}{lccccc}
Codifferential&$H^0$&$H^1$&$H^2$&$H^3$&$H^4$\\ \hline \\
$d_1=\psa{ef}e-\psa{fe}f+\psa{ee}f-\psa{ff}f$&1&0&0&0&0\\
$d_2=\psa{ff}f$&2&1&1&1&1\\
$d_3=\psa{ef}e-\psa{ff}f$&0&0&0&0&0\\
$d_4=\psa{fe}e+\psa{ff}f$&0&0&0&0&0\\
$d_5=\psa{ef}e-\psa{fe}e-\psa{ff}f$&2&2&2&2&2\\
$d_6=\psa{ee}f$&1&1&2&2&1\\\\ \hline
\end{tabular}
\end{center}
\label{coho table}
\caption{Cohomology of the six codifferentials on a $1|1$-dimensional space}
\end{table}

\section{Elements of the moduli space}
In this section we give a complete description of both the cohomology and the
multiplication structure generated by each codifferential. For a complete proof
of the cohomological structure see the next section. Let us suppose that
$V=\langle x,\theta\rangle$, where $x$ is even and $\theta$ is odd, and that
$W=\Pi V=\langle e,f\rangle$, where $\pi(x)=f$ and $\pi(\theta)=e$. Let
$m=\pi\inv\circ d\circ (\pi\tns\pi)$. Then $m$ is an associative algebra
structure on $V$, corresponding to the codifferential $d$. For each of the
codifferentials, we give the multiplication structure $m$ on $V$.
\begin{equation*}
\begin{array}{lllll}
d_1&x^2=-x&x\theta=-\theta&\theta x=-\theta&\theta^2=-x\\
d_2&x^2=x&x\theta=0&\theta x=0&\theta^2=0\\
d_3&x^2=-x&x\theta=0&\theta x=-\theta&\theta^2=0\\
d_4&x^2=x&x\theta=\theta&\theta x=0&\theta^2=0\\
d_5&x^2=-x&x\theta=-\theta&\theta x=-\theta&\theta^2=0\\
d_6&x^2=0&x\theta=0&\theta x=0&\theta^2=-x\\
\end{array}
\end{equation*}
Of these algebras, $d_1$, $d_2$, $d_5$ and $d_6$ are commutative, and $d_1$ and
$d_5$ are unital, with unit $1=-x$. In the algebras $d_2$, $d_3$, $d_4$, and $d_5$,
$\theta$ generates a nontrivial proper graded ideal, while $x$ generates a nontrivial
proper graded ideal in $d_2$ and $d_6$. Thus $d_1$ is the only simple algebra in the
moduli space. The algebras $d_3$ and $d_4$ are non-commutative and non-unital.

The algebras $d_2$, $d_3$, $d_4$ and $d_5$ are all extensions of the simple $0|1$-dimensional
associative algebra (whose structure is just the associative algebra structure of $C$). In
fact, they fit a certain pattern of extensions. The algebras $d_3$ and $d_4$ are opposite
algebras, and they are rigid in the cohomological sense. These two rigid algebras are just
the first in a sequence of rigid extensions of the $0|1$-dimensional simple algebra.

The algebra $d_5$ is the unique
extension of the simple $0|1$-dimensional algebra by the trivial $1|0$-dimensional
algebra as a unital algebra. The algebra $d_2$ is just the direct sum of the
trivial $1|0$-dimensional algebra and the simple $0|1$-dimensional algebra.

Finally, the algebra $d_6$ is an extension of the trivial $1|0$-dimensional algebra by
the trivial $0|1$-dimensional algebra, and as a consequence, it is a nilpotent algebra.
By nilpotent algebra, we mean an algebra such that a power of the algebra vanishes,
which in the finite dimensional case is equivalent to the fact that every element in
this algebra is nilpotent.

We did not use the method of extensions in calculating the nonequivalent codifferentials.
In this simple case, it is easy to solve the codifferential property $[d,d]=0$, which
gives a system of quadratic coefficients, and study the action of the group of linear
automorphisms of the underlying vector space, to arrive at the six codifferentials.
However, calculating this space by extensions reveals more of its properties, and also
gives a natural manner of organizing the codifferentials.

The remainder of this paper will be concerned with calculating the cohomology of the
codifferentials. It turned out that this aspect was the most difficult, especially for
the codifferential $d_6$, the nilpotent one.
\section{Calculating the cohomology}
The cohomology of the codifferentials is given in Table \ref{coho table} above.
With the exception of $d_6$, the pattern of cohomology is easily deduced from
the information in the table.  The pattern for $d_6$ is that $h^k=1$ if
$k=0,1\mod 4$ and $h^k=2$ otherwise.

For later use, we define  the following operator on $C(W)$. If $I$ is a
multi-index with $i_k\in\{e,f\}$, with $\ell(I)=m$, then define $\lambda^I:C^k\ra
C^{k+m}$ by $\lambda^I\pha Jj=\pha{IJ}j$. Note that the parity of $\lambda^I$
is the same as the parity of $I$. We abbreviate $\lambda^{\{e\}}=\lambda^e$.

We give a computation of the cohomology of the codifferentials on a case by
case basis.

The following theorem is a well known result.
\begin{thm}
 Suppose that a coboundary operator $D:C^n\ra C^{n+1}$ decomposes as $D=D'+D''$, given by the following diagram
 $$\xymatrix@C=13pt@R=15pt{
 \Cn{n}{a} \ar[d]^{D'} \ar[rd]^{D''} \\
 \Cn{n+1}{a} \ar[d]^{D'} \ar[rd]^{D''} \ar @{} [r]|{\oplus} & \Cn{n+1}{a+1} \ar[d]^{D'} \ar[rd]^{D''} \\
 \Cn{n+2}{a} \ar[d]^{D'} \ar[rd]^{D''} \ar @{} [r]|{\oplus} & \Cn{n+2}{a+1} \ar[d]^{D'} \ar[rd]^{D''}\ar@{}[r]|{\oplus} & \Cn{n+2}{a+2} \ar[d]^{D'} \ar[rd]^{D''} \\
 &&&&
  }
  $$
  for $a\le k\le n$ where $\Cn{n}{} =\Cn{n}1\oplus\dots\oplus\Cn{n}a\oplus\dots\oplus\Cn{n}{n}$, such that $D''$ is injective when $k=a$, and  $H(D'')=0$. Then $H(D)=0$ on the subcomplex
$\Cn{n}{k}$ for $k\ge a$.
\label{thm_arg}
\end{thm}
\begin{proof}
 Since $D^2=0$ we have
 \begin{align*}
  (D')^2=&0 & D'D''=&-D''D' & (D'')^2=&0
 \end{align*}
 Let $\ph\in\Cn{n}{}$ such that $\ph\in\ker D$. We write $\ph=\ph_a+\dots+\ph_n$ and obtain the relations,
 \begin{align*}
  D''(\ph_{k-1})+D'(\ph_k)=0
 \end{align*}
 The first relation we check is $D''(\ph_n)=0$; however, $H(D'')=0$ thus we can write $\ph_n=D''(\alpha_{n-1})$ for some $\alpha_{n-1}\in\Cn{n-1}{n-1}$. Assume we have shown
 \begin{align*}
  \ph_{k+1}=D'(\alpha_{k+1})+D''(\alpha_k).
 \end{align*}
 Then
 \begin{align*}
  0=& D''(\ph_k)+D'(\ph_{k+1}) \\
   =& D''(\ph_k)+D'D''(\alpha_k) \\
   =& D''(\ph_k-D'(\alpha_k)) \\
 \end{align*}
 Therefore we can write $\ph_k=D''(\alpha_{k-1})+D'(\alpha_{k})$ for some $\alpha_{k-1}\in\Cn{n-1}{k-1}$. This works until $k=a$. But we do know $\ph_{a+1}=D'(\alpha_{a+1})+D''(\alpha_a)$ so
 \begin{align*}
  0=&D''(\ph_a)+D'(\ph_{a+1}) \\
   =&D''(\ph_a)+D'D''(\alpha_a) \\
   =&D''(\ph_a-D'(\alpha_a))
 \end{align*}
 But $D''$ is injective when $k=a$ so $\ph_a=D'(\alpha_a)$.

\end{proof}

\subsection{$d_1=\psi^{ef}_e-\psi^{fe}_e+\psi^{ee}_f-\psi^{ff}_f$}\label{11case5}
We begin by computing the bracket of $d$ with a general element in $E^n$ and
$F^n$.
\begin{align*}
D(\ph^I_e)=&
(-1)^{I+1}\ph^I_e\psi^{ee}_f+\ph^{If}_e+(-1)^{I+1}\ph^{fI}_e+(-1)^{I+1}\ph^I_e\psi^{ef}_e+\\
&+(-1)^I\ph^I_e\psi^{fe}_e+(-1)^I\ph^I_e\psi^{ff}_f+\ph^{Ie}_f+\ph^{eI}_f\\
D(\ph^I_f)=&(-1)^I\ph^I_f\psi^{ee}_f-\ph^{If}_f+(-1)^I\ph^{fI}_f+(-1)^I\ph^I_f\psi^{ef}_e+\\
&+(-1)^{I+1}\ph^I_f\psi^{fe}_e+(-1)^{I+1}\ph^I_f\psi^{ff}_f+\ph^{eI}_e-\ph^{Ie}_e
\end{align*}

For this case, we need a different definition of $E^n$ and $F^n$.
\begin{align*}
E^n_k=&\langle\ph^I_e|{\rm The\ number\ of\ f's\ in\ I\ is\ k}\rangle\\
F^n_k=&\langle\ph^I_f|{\rm The\ number\ of\ f's\ in\ I\ is\ k}\rangle.
\end{align*}
We decompose $D$ on $E^n_k$ and $F^n_k$ as follows:
\begin{align*}
 D=&
 D'_e+D''_e+D_f:E^n_k\ra E^{n+1}_{k-1}\oplus E^{n+1}_{k+1} \oplus F^{n+1}_k\\
 D=&
 D'_2+D''_2+D_1:F^n_k\ra F^{n+1}_{k-1}\oplus F^{n+1}_{k+1} \oplus E^{n+1}_k.
\end{align*}

Since $D^2=0$, we obtain the following relations:
$$\begin{array}{ll}
(D'_e)^2=0  &(D'_2)^2=0 \\
(D''_e)^2=0  &(D''_2)^2=0 \\
D_fD'_e=-D'_2D_f & D_1D'_2=-D'_eD_1 \\
D_fD''_e=-D''_2D_f & D_1D''_2=-D''_eD_1 \\
D''_eD'_e+D'_eD''_e+D_1D_f=0 & D''_2D'_2+D'_2D''_2+D_fD_1=0
\end{array}$$

Let $\ph\in E^n$ and $\xi\in F^n$ be such that $\ph+\xi\in\ker(D)$, then we can
write $\ph=\ph_0+\ph_1+\dots+\ph_n$ and $\xi=\xi_0+\xi_1+\dots+\xi_n$. For
simplicity of notation, set $\ph_i=0$ and $\xi=0$ if $i$ is not between $1$ and
$n$. Then, for $k=0\dots n+1$, we have
\begin{align*}
D'_e(\ph_{k+1})+D''_e(\ph_{k-1})+D_1(\xi_k)=0,\qquad
D'_2(\xi_{k+1})+D''_2(\xi_{k-1})+D_f(\ph_k)=0
\end{align*}
We would like to show that $\ph+\xi=D(\eta+\alpha)$ for some $\eta\in E^{n-1}$
and $\alpha\in F^{n-1}$. This happens if, for $k=1\dots n$,
\begin{align*}
D'_e(\eta_{k+1})+D''_e(\eta_{k-1})+D_1(\alpha_k)=\ph_k,\qquad
D'_2(\alpha_{k+1})+D''_2(\alpha_{k-1})+D_f(\eta_k)=\xi_k
\end{align*}
First, note that for $k=n+1$, the equations on $\ph$ and $\xi$ reduce to
\begin{equation*}
D''_e(\ph_{n})=0,\qquad D''_2(\xi_n)=0.
\end{equation*}
Note that the coboundary operators $D''_e$ and $D''_2$ have already been
studied in the previous case, and it was shown that $H^k(D''_e)=H^k(D''_2)=0$
if $k\ge 1$. Let us suppose that we have shown that
\begin{align*}
\ph_{k+1}=&D'_e(\eta_{k+2})+D''_e(\eta_k)+D_1(\alpha_{k+1}) \\
\xi_{k+1}=&D'_2(\alpha_{k+2})+D''_2(\alpha_k)+D_f(\eta_{k+1}) \\
\ph_{k+2}=&D'_e(\eta_{k+3})+D''_e(\eta_{k+1})+D_1(\alpha_{k+2}) \\
\xi_{k+2}=&D'_2(\alpha_{k+3})+D''_2(\alpha_{k+1})+D_f(\eta_{k+2})
\end{align*}
These formulas are trivial if $k>n$. We show that if $k>0$, then we can
construct $\eta_{k-1}$ and $\alpha_{k-1}$ so that the corresponding formulas
hold for $\ph_k$ and $\xi_k$. But
$$\begin{array}{rcl}
0&=&D'_e(\ph_{k+2})+D''_e(\ph_k)+D_1(\xi_{k+1})\\
&=&D'_eD''_e(\eta_{k+1})+D'_eD_1(\alpha_{k+2})+D''_e(\ph_k)+\\
&&D_1D^{\prime}_2(\alpha_{k+2})+D_1D''_2(\alpha_k)+D_1D_f(\eta_{k+1})\\
&=&(D'_eD''_e+D_1D_f)(\eta_{k+1})-D_1D'_2(\alpha_{k+2})+\\
&&D''_e(\ph_k)+D_1D'_2(\alpha_{k+2})-D''_eD_1(\alpha_k)\\
&=&D''_eD'_e(\eta_{k+1})+D''_e(\ph_k)-D''_eD_1(\alpha_k)\\
&=&D''_e(\ph_k-D'_e(\eta_{k+1})-D_1(\alpha_k))
\end{array}$$
Thus there is a $\eta_{k-1}$ such that
$\ph_k-D'_e(\eta_{k+1})-D_1(\alpha_k)=D''(\eta_{k-1})$. A similar argument
holds for $\xi_k.$ The argument holds as long as $k>1$. When $k=1$, we can use
the same argument to show that $D''_e(\ph_1-D'_e(\eta_2)-D_1(\alpha_1))=0$, but
here, what this implies is that $\ph_1-D'_e(\eta_2)-D_1(\alpha_1)=0$, since
$D''_e$ is injective on $E^n_1$. Since $D''(\eta_0)=0$ for any $\eta_0$, this
condition is independent of the choice of $\eta_0$.

Similarly, $\xi_1=D_2'(\alpha_2)+D_f(\eta_1)$, so
\begin{align*}
0=&D_2'(\xi_1)+D_f(\ph_0)\\
=&D_2'D_f(\eta_1)+D_f(\ph_0)\\
=&
-D_fD_e'(\eta_1)+D_f(\ph_0)\\
=&D_f(\ph_0-D_e'(\eta_1)).
\end{align*}
Since $D_f$ is injective, it follows that $\ph_0=D_e'(\eta_1)$. Since $\dim
F_0^n=\dim E^{n-1}_0$ if $n>0$, $D_f$ is surjective, and we have
$\xi_0-D_2'(\alpha_1)=D_f(\eta_0)$ for some $\eta_0$. Thus when $n>0$, we have
shown that $\ph+\xi=D(\eta+\alpha)$. When $n=0$,  $E^{n-1}_0=0$, so $D_f$ is no
longer surjective. It follows that
$$H^n=\left\{
\begin{array}{ll}
\ph_f, & n=0; \\
0, & \text{otherwise}.
\end{array}
\right.$$

\subsection{$d_2=\psi^{ff}_f$}\label{11case1}
We begin by computing the coboundary of representatives from the $E^n$ and
$F^n$ space:
\begin{align*}
D(\ph^I_e)=&(-1)^{I+1}\ph^I_e\psi^{ff}_f\\
D(\ph^I_f)=&\ph^{If}_f+(-1)^{I+1}\ph^{fI}_f+(-1)^I\ph^I_f\psi^{ff}_f
\end{align*}
For $k=0\dots n-1$, let $E^n_k=\langle\ph^{f^keI}_e|\ell(I)=n-k-1\rangle$, and
let $P^n=\langle\ph^{f^n}_e\rangle.$ Then $E^n=P^n\oplus E^n_0\oplus
E^n_1\oplus \dots \oplus E^n_{n-1}$. From the formulas for $D$ above, we see
that $D:E^n\ra E^{n+1}$. More specifically, we have
 $$D=D^\prime+D^{\prime\prime}:E^n_k\ra E^{n+1}_k\oplus
 E^{n+1}_{k+1}.$$
In fact,
\begin{equation*}
D(\ph^{f^keI}_e)=(-1)^{I+1}\lambda^{f^k}D(\ph^{eI}_e)+(-1)^{I+k+1} \left\{
\begin{array}{ll}
\ph^{f^{k+1}eI}_e, & \hbox{k odd;} \\
0, & \hbox{k even.}
\end{array}
\right.
\end{equation*}
so that in particular,
\begin{equation*}
D''(\ph^{f^keI}_e)=\left\{
\begin{array}{ll}
(-1)^I \ph^{f^{k+1}eI}_e, & \hbox{k odd;} \\
0, & \hbox{k even.}
\end{array}
\right.
\end{equation*}
We obtain the following relations on $D'$ and $D''$:
\begin{equation*}
(D')^2 =0,\quad
D'D'' = - D''D',\quad
(D'')^2 =0
\end{equation*}
We also have $D:P^n\ra P^{n+1}$, and
\begin{equation*}
D(\ph^{f^n}_e)= \left\{
\begin{array}{ll}
0, & \hbox{n even;} \\
\ph^{f^{n+1}}_e, & \hbox{n odd.}
\end{array}
\right.
\end{equation*}
Note that $D''$ vanishes on $E^n_0$. Using this information we obtain the
following diagram

\begin{figure}[h]
$$\xymatrix@C=13pt@R=15pt{
 {P^0} \ar[d]^{D} \\
 {P^1 \ar @{} [r]|{\oplus}} \ar[d]^{D} &{E^1_0} \ar[d]^{D} \\
 {P^2 \ar @{} [r]|{\oplus}} \ar[d]^{D} &{E^2_0 \ar @{} [r]|{\oplus}} \ar[d]^{D} & {E^2_1} \ar[d]^{D^\prime} \ar[dr]^{D^{\prime\prime}}\\
 {P^3 \ar @{} [r]|{\oplus}} \ar[d]^{D} &{E^3_0 \ar @{} [r]|{\oplus}} \ar[d]^{D} & {E^3_1 \ar @{} [r]|{\oplus}} \ar[d]^{D^\prime} \ar[dr]^{D^{\prime\prime}} & {E^3_2} \ar[d]^{D^\prime} \ar[dr]^{D^{\prime\prime}} \\
{P^4 \ar @{} [r]|{\oplus} \ar @{} [d]|{\vdots}} & {E^4_0\ar @{} [r]|{\oplus} \ar @{} [d]|{\vdots}} & {E^4_1\ar @{} [r]|{\oplus} \ar @{} [d]|{\vdots}} &{E^4_2\ar @{} [r]|{\oplus} \ar @{} [d]|{\vdots}} &{ E^4_3 \ar @{} [d]|{\vdots}} \ar @{} [rd]|{\ddots}\\
&&&&&
}$$
\caption{Decomposition of the action of $D$ on $E$}
\label{Diagram 1}
\end{figure}

From our calculations on the $P$ space, we see that $D$ oscillates between the
zero map and an isomorphism. This means that when $n\geq 1$, the cohomology on
the $P$ space vanishes, but when $n=0$ we have $H^0(P)=\langle\ph_e\rangle$.
\begin{rem}
Given a complex $C_k$ with coboundary operators $D:C_k\ra C_{k+1}$, if the maps alternate between the zero
map and isomorphisms, then if the initial map is an isomorphism, the cohomology will vanish, and if the initial
map is the zero map, then only cohomology is the initial space.
 \label{rem_arg}
\end{rem}

By \thmref{thm_arg} with $a=1$, the cohomology vanishes on $E^n_k$ for $k\ge1$. Finally, we compute the cohomology on the complex $E^n_0$. First, note that $\lambda^e:E^n\mapsiso \tilde{E}^{n+1}_0$ is an isomorphism commuting
with $D$, which implies that we have a map $\lambda^e_*:H^n(E)\mapsiso H^{n+1}(E_0).$ For $n\geq2$ we have $H^n(E)= H^n(E_0)\cong H^{n-1}(E)$, so that $\lambda^e_*:H^{n-1}(E)\mapsiso H^n(E)$. Moreover, $\lambda^e$ maps $P^0$ isomorphically to $E^1_0$, and since $D$ is trivial on $P^0$ and $E^1_0$, it also commutes with $D$. Thus $H^0(E)=H^0(P)=H^1(E^1_0)$. It follows that $H^n(E)=\langle(\lambda^e)^n\ph_e\rangle=\langle\ph^{e^n}_e\rangle.$

Next, we consider the cohomology of $D$ on the $F$ space. As before we have
$F^n=Q^n\oplus F^n_0\oplus \dots\oplus F^n_{n-1}$ where
 $$F^n_k=\langle\ph^{f^keI}_f|\ell(I)=n-k-1\rangle{\rm\quad
 and\quad}Q^n=\langle\ph^{f^n}_f\rangle$$
Now we redefine our maps, and decompose $D$ as $D=D^\prime+D^{\prime\prime}:F^n_k\ra F^{n+1}_k\oplus
F^{n+1}_{k+1}$, where
$$D''(\ph^{f^keI}_f)=\left\{
\begin{array}{ll}
(-1)^{I+1}\ph^{f^{k+1}eI}_f, & \hbox{k even;} \\
0, & \hbox{k odd.}
\end{array}
\right.$$
We see the cohomology of the $D''$ map on $F$ vanishes completely. We obtain
the following diagram:

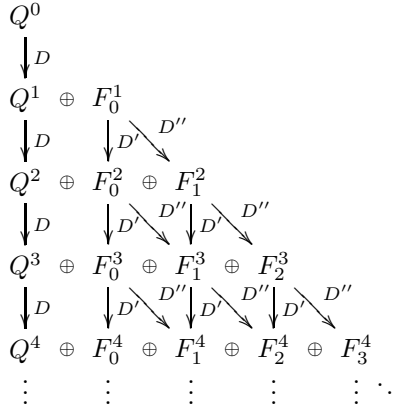
\begin{figure}[ht]
$$\xymatrix@C=13pt@R=15pt{
 {Q^0} \ar[d]^{D} \\
 {Q^1 \ar @{} [r]|{\oplus}} \ar[d]^{D} &{F^1_0} \ar[d]^{D^\prime} \ar[dr]^{D^{\prime\prime}} \\
 {Q^2 \ar @{} [r]|{\oplus}} \ar[d]^{D} &{F^2_0 \ar @{} [r]|{\oplus}} \ar[d]^{D^\prime} \ar[rd]^{D^{\prime\prime}} & {F^2_1} \ar[d]^{D^\prime} \ar[dr]^{D^{\prime\prime}}\\
 {Q^3 \ar @{} [r]|{\oplus}} \ar[d]^{D} &{F^3_0 \ar @{} [r]|{\oplus}} \ar[d]^{D^\prime} \ar[rd]^{D^{\prime\prime}} & {F^3_1 \ar @{} [r]|{\oplus}} \ar[d]^{D^\prime} \ar[dr]^{D^{\prime\prime}} & {F^3_2} \ar[d]^{D^\prime} \ar[dr]^{D^{\prime\prime}} \\
{Q^4 \ar @{} [r]|{\oplus} \ar @{} [d]|{\vdots}} & {F^4_0\ar @{} [r]|{\oplus} \ar @{} [d]|{\vdots}} & {F^4_1\ar @{} [r]|{\oplus} \ar @{} [d]|{\vdots}} &{F^4_2\ar @{} [r]|{\oplus} \ar @{} [d]|{\vdots}} &{ F^4_3 \ar @{} [d]|{\vdots}} \ar @{} [rd]|{\ddots}\\
&&&&& }$$
\caption{Decomposition of the action of $D$ on $F$} \label{Diagram 2}
\end{figure}

Using \thmref{thm_arg} we see that the cohomology on the $F^n_k$ spaces vanishes completely.

The only space we have left is the $Q^n=\langle\ph^{f^n}_f\rangle$.
 $$D(\ph^{f^n}_f)=\left\{
                    \begin{array}{ll}
                      \ph^{f^{n+1}}_f, & \hbox{n odd;} \\
                      0, & \hbox{n even.}
                    \end{array}
                  \right.$$
We see the map alternates between the zero map and an isomorphism, meaning
the cohomology on the $Q^n$ space is zero for $n\geq 1$ and
$H^0(Q)=\langle\psi_f\rangle$.

Thus we see the total cohomology is given by
 $$H^n(d_1)=\left\{
         \begin{array}{ll}
           \langle\phi_e,\psi_f\rangle, & n=0; \\
           \langle\phi^{e^n}_e\rangle, & n\geq 1.\\
         \end{array}
       \right.$$

\subsection{$d_3=\psi^{ef}_e-\psi^{ff}_f$}\label{11case2}
We begin by computing the bracket of $d$ with a general element in $E^n$ and
$F^n$.
\begin{align*}
D(\ph^I_e)=&\ph^{If}_e+(-1)^{I+1}\ph^I_e\psi^{ef}_e+(-1)^I\ph^I_e\psi^{ff}_f\\
D(\ph^I_f)=&\ph^{eI}_e-\psi^{If}_f+(-1)^I
\psi^{fI}_f+(-1)^I\ph^I_f\psi^{ef}_e+(-1)^{I+1}\ph^I_f\psi^{ff}_f
\end{align*}
We see that
\begin{align*}
D:&E^n\ra E^{n+1}\\
D=D_e+D_f:&F^n\ra E^{n+1}\oplus F^{n+1}.
\end{align*}
Since $D^2=0$, we have the following relations:
 $$DD_e=-D_eD_f,\qquad D^2_f=0.$$
Note that $D_e(\ph^I_f)=\ph^{eI}_e$, so that $\im(D_e) = E^n_0$ this will be considered later.
We define $E^n_k$ analogous to our previous case, so we decompose $E^n=P^n\oplus E^n_0\oplus\dots\oplus
E^n_{n-1}$. Being more specific on the $E$ space we see
\begin{align*}
 D=D'+D'':E^n_k&\ra E^{n+1}_k\oplus E^{n+1}_{k+1} \\
 D:P^n&\ra P^{n+1}
\end{align*}
and $D$: $E^n_0\ra E^{n+1}_0$. Thus, by \thmref{thm_arg} the cohomology on $E^n_k$ for $k\ge1$, and by calculation and Remark \ref{rem_arg} the cohomology vanishes on $P^n$.

Let $\ph\in E^n$ and $\xi\in F^n$ such that $\ph+\xi\in\ker(D)$. Then the
following must hold
 $$D(\ph)=-D_e(\xi), \quad \text{and} \quad D_f(\xi)=0.$$

However we see that the second equality follows directly from the first,
\begin{align*}
D_eD_f(\xi) &=-DD_e(\xi) \\
            &=DD(\ph) \\
            &=0.
\end{align*}
Since $D_e$ is injective it follows that $D_f(\xi)=0.$

Because the cohomology vanishes on $P^n$ and $E^n_k$ $k\geq 1$, we will now
consider $\ph$ to be an element of $E^n_0$ and show that it is a coboundary.
Since $D_e$ is an isomorphism we have $\ph=D_e(\xi^\prime)$, for some
$\xi^\prime\in F^{n-1}$ and it follows that
\begin{align*}
D_e(\xi) &=-D(\ph)\\
         &=-DD_e(\xi^\prime)\\
         &=D_eD_f(\xi^\prime)
\end{align*}
So $\xi=D_f(\xi^\prime)$ thus $\ph+\xi=D(\xi^\prime)$. Since $\ph+\xi$ is in
the kernel and $\ph+\xi\in \im(D)$ then the cohomology vanishes when $n>1$.
This argument breaks down in $E^1_0$. However, since
$D(\ph^e_e)=2\ph^{ef}_e\ne0$, it follows that there are no cocycles in $E^n_0$.
 $$H^n(d_2)=0\qquad\text{for all } n$$

\subsection{$d_4=\psi^{fe}_e+\psi^{ff}_f$}\label{11case3}
This case is completely analogous to case 2, if we define
$E^n_k=\langle\ph^{Ief^k}_e|\ell(I)=n-k-1\rangle$ and similarly for $F^n_k$.

\subsection{$d_5=\psi^{ef}_e-\psi^{fe}_e-\psi^{ff}$}\label{11case4}
In this case, we have
\begin{align*}
D(\ph^I_e)=&
 \ph^{If}_e+(-1)^{I+1}\ph^{fI}_e+(-1)^I\ph^I_e\psi^{ff}_f+(-1)^{I+1}
 \ph^I_e\psi^{ef}_e+(-1)^I\ph^I_e\psi^{fe}_e\\
D(\ph^I_f)=&-\ph^{If}_f+(-1)^I\ph^{fI}_f+\ph^{eI}_e-\ph^{Ie}_e
 \\&+(-1)^{I+1}\ph^I_f\psi^{ff}_f+(-1)^I\ph^I_f\psi^{ef}_e+(-1)^{I+1}\ph^I_f\psi^{fe}_e.
\end{align*}
Thus
\begin{align*}
&D:E^n\ra E^{n+1}\\
&D=D_e+D_f:F^n\ra E^{n+1}\oplus F^{n+1}
\end{align*}
For this case, we need to split $E^n_0$ and $F^n_0$ into smaller pieces. Define
\begin{align*}
S^n=&\langle\ph^{e^n}_e\rangle,\tilde{E}^n_0=\{\ph^{eI}_e|\ell(I)=n-1,I\ne e^{n-1}\}\\
T^n=&\langle\ph^{e^n}_f\rangle, \tilde{F}^n_0=\{\ph^{eI}_f|\ell(I)=n-1,I\ne e^{n-1}\},
\end{align*}
so that $E^n_0=\tilde{E}^n_0\oplus S^n$ and $F^n_0=\tilde{F}^n_0\oplus T^n$.
Let
\begin{align*}
\tilde{E}^n=&P^n\oplus \tilde{E}_0^n\oplus E^n_1\oplus\dots\oplus E^n_{n-1}\\
\tilde{F}^n=&Q^n\oplus \tilde{F}_0^n\oplus F^n_1\oplus\dots\oplus F^n_{n-1},
\end{align*}
so that $E^n=S^n\oplus \tilde{E}^n$ and $F^n=T^n\oplus\tilde{F}^n$

We decompose the action of $D$ on $E^n$ and $F^n$ as follows:
\begin{align*}
D=D'+D'':& E^n_k\ra E^{n+1}_k\oplus E^{n+1}_{k+1}\\
D=D'_f+D''_f+D^0_e+D'_e:& F^n_k\ra F^{n+1}_k \oplus F^{n+1}_{k+1}\oplus E^{n+1}_0\oplus E^{n+1}_k.
\end{align*}
When $k=0$, there is some ambiguity about the maps $D^0_e$ and $D'_e$, which we
resolve by taking $D_e=D'_e$ on $E^n_0$.

Note that
\begin{equation*}
D''(\pha{f^keI}e)=\left\{
\begin{array}{ll}
(-1)^{I+1}\ph^{f^{k+1}eI}_e, &k\text{ odd;} \\
0, & k \text{ even.}
\end{array}
\right.
\end{equation*}
From this we see that $D''D^0_e=0$. Since $D^2=0$, we have the following
relations
\begin{align*}
(D')^2=&0,\qquad &(D'')^2=&0,\qquad &D'D''=-D''D'\\
(D'_f)^2=&0,\qquad &(D''_f)^2=&0,\qquad &D'_fD''_f=-D_f''D_f'\\
D^0_eD_f=&-D'D^0_e,\qquad &D'_eD''_f=&-D''D'_e,\qquad &D'_eD'_f=-D'D'_e.
\end{align*}

We will show the cohomology $H_f$ of $D_f$ on the $\tilde{F}$ space vanishes.
First, we determine what $D_f$ does to the $Q^n$ space.
\begin{equation*}
D_f(\ph^{f^n}_f)=\left\{
\begin{array}{ll}
0, & \hbox{n even;} \\
-\ph^{f^n}_f, & \hbox{n odd.}
\end{array}
\right.
\end{equation*}
Thus $D_f:Q^n\ra Q^{n+1}$. Furthermore we see that the cohomology vanishes on
the subcomplex $Q^n$ for $n\geq1$.

Next, note that
\begin{equation*}
D''_f(\ph^{f^keI}_f)=\left\{
\begin{array}{ll}
0, & \hbox{ k even;} \\
\ph^{f^{k+1}eI}_f, & \hbox{k odd}
\end{array}
\right.
\end{equation*}
Thus we have the following diagram
$$\xymatrix@C=13pt@R=15pt{
 {T^0} \ar[d]^{0} \\
 {T^1 \ar @{} [r]|{\oplus}} \ar[d]^{0} &{Q^1} \ar[d]^{D} \\
 {T^2 \ar @{} [r]|{\oplus}} \ar[d]^{0} &{Q^2 \ar @{} [r]|{\oplus}} \ar[d]^{D}&{\tilde{F}^2_0 \ar @{} [r]|{\oplus}} \ar[d]^{D} & {F^2_1} \ar[d]^{D'} \ar[dr]^{D''}\\
 {T^3 \ar @{} [r]|{\oplus}} \ar[d]^{0} &{Q^3 \ar @{} [r]|{\oplus}} \ar[d]^{D} &{\tilde{F}^3_0 \ar @{} [r]|{\oplus}} \ar[d]^{D}& {F^3_1 \ar @{} [r]|{\oplus}} \ar[d]^{D'} \ar[dr]^{D''} & {F^3_2} \ar[d]^{D'} \ar[dr]^{D''} \\
 {T^4 \ar @{} [r]|{\oplus} \ar @{} [d]|{\vdots}} & {Q^4\ar @{} [r]|{\oplus} \ar @{} [d]|{\vdots}}&{\tilde{F}^4_0 \ar @{} [r]|{\oplus} \ar @{} [d]|{\vdots}} & {F^4_1\ar @{} [r]|{\oplus} \ar @{} [d]|{\vdots}} &{F^4_2\ar @{} [r]|{\oplus} \ar @{} [d]|{\vdots}} &{ F^4_3 \ar @{} [d]|{\vdots}} \ar @{} [rd]|{\ddots}\\
 &&&&&&
}$$

Now we show that $\lambda^e$ commutes with the $D_f$ operator:
\begin{align*}
D_f(\lambda^e\ph^I_f)=&
-\lambda^e\ph^{If}_f+(-1)^I\ph^{feI}_f+
(-1)^{I+1}\lambda^e\ph^I_f\psi^{ff}_f+(-1)^I\lambda^e\ph^I_f\psi^{ef}_e\\
&+(-1)^I\lambda^e\ph^{fI}_f+(-1)^{I+1}\lambda^e\ph^I_f\psi^{fe}_e+(-1)^{I+1}\ph^{feI}_f\\
=&\lambda^eD_f(\ph^I_f).
\end{align*}

Thus $D_f:\tilde{F}^n_0\ra \tilde{F}^{n+1}_0.$ We also note that
$\lambda^e:\tilde{F}^n\mapsiso \tilde{F}^{n+1}_0$ for $n\geq1$, and since
$\lambda^e$ commutes with $D_f$, we have an isomorphism
$\lambda^e_\star:H^n_f(\tilde{F})\mapsiso H^{n+1}_f(\tilde{F}_0)$ for $n\geq1$.
Now $\tilde{F}^1=Q^1$, and $D_f$ is an isomorphism $Q^1\ra Q^2$, so
$H_f^1(\tilde{F})=0$. Thus $H^n_f(\tilde F)=0$.

By a similar argument, we see that the cohomology $H(\tilde E)$ induced by $D$
acting on the subcomplex $\tilde{E}$ also vanishes. In fact, one shows that
$\lambda^e$ commutes with $D$ on $\tilde{E}$, and therefore $H^n(\tilde
E)=H^{n+1}(\tilde{E}_0)$, just as in the $\tilde F$ space.

Now consider the action of $D$ on the space $\tilde E\oplus\tilde F$. Suppose
$\ph\in \tilde{E}$ and $\xi\in\tilde{F}$ are such that $D(\ph+\xi)=0$. Then
$D(\ph)+D_e(\xi)=0$ and $D_f(\xi)=0$. In fact, the second equality follows from
the first, because $D_e$ is injective on $\tilde F$, so if the first equality
holds, then $D_eD_f(\xi)=-DD_e(\xi)-D^2(\ph)=0$. (Here we used the relation
$D_eD_f=-DD_e$, which follows immediately from the fact that $D^2=0$.) Since
$D_e$ is injective, it follows that $D_f(\xi)=0$. Since $D_f$ has trivial
cohomology, we must have $\xi=D_f(\eta)$ for some $\eta\in\tilde F$. Therefore
\begin{equation*}
0=D(\ph)+D_e(\xi)=D(\ph)+D_eD_f(\eta)=D(\ph)-DD_e(\eta)=D(\ph-D_e(\eta)).
\end{equation*}
But the cohomology of $D$ on $\tilde E$ vanishes, so this means that
$\ph-D_e(\eta)=D(\tau)$ for some $\tau\in\tilde E$. But then
$D(\eta+\tau)=\ph+\xi$. This shows that the cohomology of $D$ on $\tilde
E\oplus \tilde F$ vanishes.

Finally we are left with the $T^n$ and $S^n$ spaces. The important fact is that
$D$ vanishes on $T^n$ and $S^n$:
\begin{align*}
D(\ph^{e^n}_e)&=\ph^{e^nf}_e-\ph^{fe^n}_e-\ph^{efe^{n-1}}_e-\dots-
\ph^{e^nf}_e+\ph^{fe^n}_e+\dots+\ph^{e^{n-1}fe}_e\\
&=0.
\end{align*}.
Similarly $D(\pha{e^n}f)=0$. Thus both maps are the zero maps meaning the
cohomology is exactly these two spaces. More precisely, we have
$H^n(d_4)=\langle\phi^{e^n}_e,\psi^{e^n}_f\rangle$.

\subsection{$d_6=\psi^{ee}_f$}\label{11case6}
For case six, the methods employed in the previous cases completely fail, and
as a result, a different method is necessary. Once again, we begin by computing
the bracket of \emph{$d$} with a general element in $E^n$ and $F^n$.
\begin{align*}
D(\ph^I_e)&=\ph^{eI}_f+\ph^{Ie}_f+(-1)^I\ph^I_e\ph^{ee}_f\\
D(\ph^I_f)&=(-1)^{I+1}\ph^I_f\ph^{ee}_f
\end{align*}
Therefore, we have decompositions
\begin{align*}
D=D_e+D_f&:E^n\ra E^{n+1}\oplus F^{n+1}\\
D&:F^n\ra F^{n+1}
\end{align*}
Note that $D_f$ is injective and that
\begin{equation*}
D_e^2=0,\qquad D_fD_e=-DD_f.
\end{equation*}

Thus $D_e$ is a coboundary operator on $E$, giving a cohomology
$H^n_e=H^n(D_e)$. We first compute this cohomology, and use it to compute the
cohomology in general. Define the Decleene map
$\theta=\lambda^{ef}-\lambda^{fe}$. Then we claim that $\theta$ commutes with
$D_e$ on $E$, and with $D$ on $F$. To see this, note that
\begin{align*}
D_e\theta(\pha Ie)&=D_e(\pha{efI}e-\pha{feI}e)\\
&=\s{I}\pha{efI}e\psa{ee}f+\s{I+1}\pha{feI}e\pha{ee}f\\
&=\s{I}\pha{eeeI}e+\s{I+1}\lambda^{ef}\pha Ie\psa{ee}f+
\s{I+1}\pha{eeeI}e-\s{I+1}\lambda^{ef}\pha Ie\psa{ee}f\\
&=\theta D_e(\pha Ie).
\end{align*}
The proof that the DeCleene map commutes with $D$ on $F$ is similar. In fact,
note that the action of $D$ on $F$ is essentially the same as $D_e$ on $E$.

Next, note that if $\ph$ is a $D_e$-coboundary,  then every term in $\ph$ must
have a double $e$. Therefore, any $D_e$-cocycle which has a term without a
double $e$ must be nontrivial.  In particular, the 0-cocycle $\ph_e$ and the
1-cocycle $\pha ee$ are nontrivial $D_e$-cocycles.  Define the Decleene
cocycle $\dec^n_e$ and $\dec^n_f$ by
\begin{align*}
\dec^{2n}_e&=\theta^n\ph_e, &\dec^{2n+1}_e&=\theta^n\pha ee,\\
\dec^{2n}_f&=\theta^n\ph_f, &\dec^{2n+1}_f&=\theta^n\pha ef.
\end{align*}
Then $\dec^n_e$ is a nontrivial $D_e$-cocycle. Also $\dec^n_f$ is nontrivial if
we consider only the cohomology of $D$ restricted only to the $F$ space. We
shall discuss later when it is a nontrivial cocycle on the whole space
$C=E\oplus F$.

Let $B^n_e$ be the space of $D_e$ $n$-coboundaries, $Z^n_e$ be the
$n$-cocycles, $z_n=\dim(Z^n_e)$, $b_n=\dim(B^n_e)$ and $h_n=\dim H^n_e$. Then
$h_n=z_n-b_n$ and $z_n+b_{n+1}=2^n$. Because there is a nontrivial Decleene
cocycle in every degree, we know that $h_n\ge 1$. We wish to show that $h_n=1$.

To see this, note that $\theta_n:B^n\ra B^{n+2}$ is injective,
$D_e\circ\lambda^f:E^n\ra B^{n+2}$ is also injective, and the images of these
operators are independent subspaces.  As a consequence, we must have
$b_{n+2}\ge b_n+2^n$. We will show that $b_n+b_{n+1}\ge 2^n-1$ for $n\ge 0$.
Since $b_0=b_1=0$, and $b_2=1$ by direct computation, the formula holds for
$n\le 1$. Suppose it holds for $k<n$, and that $n\ge 2$. Then
\begin{equation*}
b_{n}+b_{n+1}\ge b_{n-2}+2^{n-2}+b_{n-1}+2^{n-1}\ge 2^{n-2}+2^{n-2}+2^{n-1}=2^{n}.
\end{equation*}
Thus, by induction, the formula  holds for all $n$. Using this formula, we
obtain
\begin{align*}
1\le h_n=z_n-b_n=2^n-b_{n+1}-b_n\le 2^n-(2^n-1)=1.
\end{align*}

First, let us say that $\dec^n_e$ extends to a $D$-cocycle if there is some
$\eta\in F^n$ such  that $\dec^n_e+\eta$ is a $D$-cocycle. If $\dec^n_e$
extends, then let $\dec^n=\dec^n_e+\eta$ be some arbitrary extension of
$\dec^n_e$.

Suppose that $D(\ph+\xi)=0$ for some $\ph\in E^n$ and $\xi\in F^n$. Then
$D_f(\ph)+D(\xi)=0$ and $D_e(\ph)=0$. In fact, the second equation follows from
the first one. For, suppose the first equality holds. Then
\begin{equation*}
D_fD_e(\ph)=-DD_f(\ph)=D^2(\xi)=0.
\end{equation*}
Using the fact that $D_f$ is injective, we see that $D_e(\ph)=0$. Now we can
write $\ph=a\dec^n_e+D_e(\alpha)$ for some $\alpha\in E^{n-1}$, because we know
that $h_n=1$.

Note that if $\dec^n_e$ does not extend to a $D$-cocycle, then $a=0$. This is
because
\begin{align*}
0&=D_f(\ph)+D(\xi)\\&=D_f(a\dec^n_e)+D_fD_e(\alpha)+D(\xi)\\&=D_f(a\dec^n_e)-DD_f(\alpha)+D(\xi)
\\&=D_f(a\dec^n_e)-D(D_f(\alpha)-\xi)),
\end{align*}
so that if $a\ne0$ we have $D_f(\dec^n_e)=D(\eta)$, where
$\eta=\tfrac1a(D_f(\alpha)-\xi)$. Next, we claim  that
$\ph+\xi=b\dec^n_f+D(\alpha+\beta)$ for some $\beta\in F^{n-1}$. To see this,
first suppose that $a=0$. Then
\begin{align*}
0=D(\ph+\xi)=D_fD_e(\alpha)+D(\xi)=-D(D_f(\alpha)+D(\xi)=D(\xi-D_f(\alpha)).
\end{align*}
Thus $\xi-D_f(\alpha)$ is a $D$-cocycle lying in $F^n$, which means it can be
written in the form $\xi-D_f(\alpha)=b\dec^n_f+D(\beta)$ for some $\beta\in
F^{n-1}$. But this means $\ph+\xi=b\dec^n_f+D(\alpha+\beta)$, as desired.

On the other hand, if $a\ne0$ then $a\dec^n_e=a\dec^n+a\eta$, where $\eta\in
F^{n-1}$, so $\ph=a\dec^n+a\eta+D_e(\alpha)$, and then
\begin{equation*}
0=D(\ph+\xi)=D(a\eta)+D_fD_e(\alpha)+D(\xi)=D(\xi+a\eta-D_f(\alpha)).
\end{equation*}
Thus in this case, we can express $\xi+a\eta-D_f(\alpha)=b\dec^n_f+D(\beta)$,
so we obtain
\begin{equation*}
\ph+\xi=a\dec^n+a\eta+\xi=a\dec^n+b\dec^n_f+D(\alpha+\beta).
\end{equation*}
From the equation above, it follows that the dimension of $H^n$ is at most 2,
depending on whether $\dec^n_e$ extends to a $D$-cocycle and whether $\dec^n_f$
is a nontrivial cocycle.

Now we show that the non triviality of $\Ch^n_f$ is linked to the whether or
not we can extend $\Ch^{n-1}_e$.

Suppose $\dec^n_f$ is trivial, ie. $\dec^n_f=D(\ph+\xi)$ for some $\ph\in
E^{n-1}$, $\xi\in F^{n-1}$. $D_e(\ph)=0$ so $\ph=a\dec^{n-1}_e+D_e(\alpha)$ for
some $\alpha\in E^{n-1}$. If $\Ch^{n-1}_e$ extends, so
$\dec^{n-1}_e=\dec^{n-1}+\eta$ then
\begin{align*}
\dec^n_f&=D_f(\ph)+D(\xi)=aD(\eta)+D_fD_e(\alpha)+D(\xi)\\&=
D(\eta)-DD_f(\alpha)+D(\xi)=D(\xi+\eta-D_f(\alpha)).
\end{align*}
But then $\dec^n_f$ is a coboundary in the $F$ space, which is impossible. Thus
if $\dec^n_f$ is trivial, $\Ch^{n-1}_e$ does not extend to a $D$-cocycle.

On  the other hand, suppose $\dec^{n-1}_e$ does not extend to a $D$-cocycle.
Then $D_f(\dec^{n-1}_e)$ is a $D$-cocycle, lying in $F$, which is nontrivial in
terms of the $D$-cohomology restricted to $F$, so we must have
$D_f(\dec^{n-1}_e)=a\dec_f^n+D(\beta)$ for some $\beta\in F^{n-1}$, where
$a\ne0$. But then $\dec^n_f=D(\tfrac1a(\dec^{n-1}_e-\beta))$, and therefore
$\dec^n_f$ is trivial.

We now will examine the possible cases of $n \mod 4$, showing that in each
case, we can either determine that $\dec^n_e$ does not extend, or that
$\dec^n_f$ is nontrivial. From this, we will determine precisely the dimension
of $H^n$.

First suppose that $n=0\mod 4$. If $n=0$, one can check explicitly that
$\dec^0_e$ does not extend. Otherwise, consider the terms of the form
$\ph^{(ef)^{n/2}}_e$ and $\ph^{(fe)^{n/2}}_e$ appearing in $\dec^n_e$, which
have the same sign. When applying $D_e$ to these terms, we will obtain two
terms of the form $\ph^{(ef)^{n/2}e}_e$, which appear nowhere else in the
expression for $D_f(\dec^n_e)$. These terms do not have any double $e$'s, and
therefore $D_f(\dec^n_e)$ is not a $D$-coboundary of any $\eta\in F^{n-1}$.
Therefore $\dec^n_e$ does not extend. Note that if $n=2\mod 4$, the terms will
have opposite signs, so the argument above does not apply, and if $n$ is odd,
there are no such terms. Thus this argument works precisely when $n=0\mod 4$.

Also, if $n=0\mod 4$, we claim that $\dec^n_f$ is a nontrivial cocycle. For
suppose that $\dec^n_f=D(\ph+\xi)$, where $\ph\in E^{n-1}$ and $\xi\in
F^{n-1}$. Note that this implies that $D_e(\ph)=0$, so we can express
$\ph=a\dec^{n-1}+D_e(\alpha)$ for some $\alpha\in E^{n-2}$. Now $\dec^n_f$
contains  the term $\ph^{(ef)^{n/2}}_f$. This terms contain no double $e$'s,
and therefore it cannot arise in $D(\xi)$. Therefore it must appear in
$D_f(\ph)=aD_f(\dec^{n-1}_e)+D_fD_e(\alpha)$. Certainly, it doesn't appear in
$D_fD_e(\alpha)$, because every term in  that expression has double $e$'s.
Therefore, it must appear in $D_f(\dec^{n-1}_e)$. But terms in that expression
contain more $e$'s than $f$'s, while the term $\ph^{(ef)^{n/2}}_f$ has the same
number of $e$'s and $f$'s.  This shows that $\dec^n_f$ is nontrivial.

A similar argument shows that if $n=2\mod 4$, then $\dec^n_f$ is nontrivial.
Finally, if $n=3\mod 4$, then $\dec^n_f$ contains the term
$\ph^{(ef)^{(n-1)/2}e}_f$. If we express $\dec^n_f=D(\ph+\xi)$, as before,
where $\ph=a\dec^{n-1}_e+ D_e(\alpha)$, then the previous argument shows that
this term must arise from $D_f(\dec^{n-1}_e)$. There are two terms in
$\dec^{n-1}_e$, $\ph^{(ef)^{(n-1)/2}}_e$, and $\ph^{(fe)^{(n-1)/2}}_e$, which
could contribute a term of the right type. However, these terms occur in
$\dec^{n-1}$ with opposite sign, because $(n-1)/2$ is odd. Therefore, the term
$\ph^{(ef)^{(n-1)/2}e}_f$ can not appear in $D_f(\dec^{n-1}_e)$. As a
consequence $\dec^n_f$ is nontrivial.

The chart below summarizes what we have deduced. In the chart, a zero in the
$\dec^n$ column means that $\dec^n$ does not exist, while in the $\dec^n_f$
column, a zero means that $\dec^n_f$ is a trivial cocycle.
\begin{center}
\begin{tabular}{|c|c|c|c|}
\hline
$n\mod 4$&$\Ch^n$&$\Ch^n_f$&Total\\
\hline
0&0&1&1\\
1&1&0&1\\
2&1&1&2\\
3&1&1&2\\
\hline
\end{tabular}
\end{center}
When $\dec^n$ exists, it cannot happen that $a\dec^n+b\dec^n_f=D(\ph+\xi)$,
where $\ph\in E^{n-1}$ and $\xi\in F^{n-1}$ if $a\ne 0$. For this to happen, it
would be necessary that $D_e(\ph)=a\dec^n_e$. But this is impossible, because
$\dec^n_e$ is a nontrivial $D_e$-cocycle. Therefore, we obtain that the
dimension of $H^n(d_6)$ is given by the Total column in the chart above.

\section{Infinitesimal Deformations}
To compute the infinitesimal deformations we only need consider objects with odd cohomology in degree 2. This leaves only $d_5$ and $d_6$. These will be considered separately.
\subsection{$d_5$}
The odd cohomology of $d_5$ is given by $\psa{ee}{f}$, so we determine when
$$d_t=\psi^{ef}_e-\psi^{fe}_e-\psi^{ff}+t\psa{ee}{f}$$
is isomorphic to another codifferential. We see that when $t\ne0$ this is isomorphic to $d_1$.
\subsection{$d_6$}
The odd cohomology of $d_6$ is given by $\psa{ef}e-\psa{fe}e-\psa{ff}{f}$, so we determine when
$$d_t=\psa{ee}{f}+t(\psa{ef}e-\psa{fe}e-\psa{ff}{f})$$
is isomorphic to another codifferential. We see that when $t\ne0$ this is isomorphic to $d_1$.
\section{Versal Deformations}
In this case the versal deformations coincide exactly with the infinitesimal deformations. This can be seen be taking the bracket of the infinitesimal deformation with itself. If the result is zero then the deformations coincide.
\subsection{Diagram of Deformations}
\begin{figure}[ht]
 $$\xymatrix@R=5mm@C=3mm{
 \bullet^{d_2} & \bullet^{d_3} & \bullet^{d_4} & ^{d_5}\bullet \ar@{~>}[rd] & & \bullet^{d_6} \ar@{~>}[ld] \\
 &&&&\bullet_{d_1}&
 }$$
 \caption{Infinitesimal deformations of a $1|1$-dimensional vector space}\label{dia}
\end{figure}
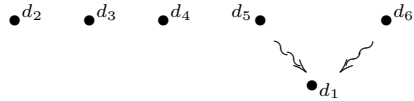

For a visual of the deformations we provide Figure \ref{dia}.

\bibliographystyle{amsplain}

\begin{thebibliography}{10}

\bibitem{bsz}
Y.~A. Bahturin, S.~K. Sehgal, and M.~V. Zaicev, \emph{Finite dimensional graded
  simple algebras}, math.RA/0512202v1, 2005.

\bibitem{bppw1}
L.~Brunshidle, M.~Penkava, M.~Phillipson, and D.~Wackwitz, \emph{Fundamental
  theorem of finite dimensional graded associative algebras}, preprint, 2009.

\bibitem{fp10}
A.~Fialowski and M.~Penkava, \emph{Formal deformations, contractions and moduli
  spaces of {L}ie algebras}, International Journal of Theoretical Physics
  (2007), to appear.

\bibitem{gers}
M.~Gerstenhaber, \emph{The cohomology structure of an associative ring}, Annals
  of Mathematics \textbf{78} (1963), 267--288.

\bibitem{gers1}
\bysame, \emph{On the deformations of ringe and algebras {I}}, Annals of
  Mathematics \textbf{79} (1964), 59--103.

\bibitem{gers2}
\bysame, \emph{On the deformations of ringe and algebras {II}}, Annals of
  Mathematics \textbf{84} (1966), 1--19.

\bibitem{gers3}
\bysame, \emph{On the deformations of ringe and algebras {III}}, Annals of
  Mathematics \textbf{88} (1968), 1--34.

\bibitem{gers4}
\bysame, \emph{On the deformations of ringe and algebras {IV}}, Annals of
  Mathematics \textbf{99} (1974), 257--276.

\bibitem{hoch}
G.~Hochschild, \emph{On the cohomology groups of an associative algebra},
  Annals of Mathematics \textbf{46} (1945), 58--67.

\bibitem{Karrer}
Guido Karrer, \emph{graded division algebras}, Mathematische Zeitschrift
  \textbf{133} (1973), 67--73.

\bibitem{knus}
M.~A. Knus, \emph{Algebras graded by a group}, Lecture Notes in Mathematics,
  vol.~92, pp.~117--133, Springer Verlag, 1969.

\bibitem{mp1}
M.~Mulase and M.~Penkava, \emph{Ribbon graphs, quadratic differentials on
  {R}iemann surfaces, and algebraic curves defined over
  $\overline{\mathbb{q}}$}, Asian Journal of Mathematics \textbf{2} (1998),
  875--920, math-ph/9811024.

\bibitem{ps1}
M.~Penkava and A.~Schwarz, \emph{On some algebraic structures arising in string
  theory}, Perspectives in Mathematical Physics (R~Penner and S.~Yau, eds.),
  International Press, 1994, pp.~219--227.

\bibitem{ps2}
\bysame, \emph{\hbox{$A_\infty$} algebras and the cohomology of moduli spaces},
  Dynkin Seminar, vol. 169, American Mathematical Society, 1995, pp.~91--107.

\bibitem{pv1}
M.~Penkava and P.~Vanhaecke, \emph{Deformation quantization of polynomial
  {P}oisson algebras}, Journal of Algebra (2000), no.~227, 365--393.

\bibitem{sta4}
J.D. Stasheff, \emph{The intrinsic bracket on the deformation complex of an
  associative algebra}, Journal of Pure and Applied Algebra \textbf{89} (1993),
  231--235.

\end{thebibliography}
\providecommand{\bysame}{\leavevmode\hbox to3em{\hrulefill}\thinspace}
\providecommand{\MR}{\relax\ifhmode\unskip\space\fi MR }
\providecommand{\MRhref}[2]{%
  \href{http://www.ams.org/mathscinet-getitem?mr=#1}{#2}
}
\providecommand{\href}[2]{#2}

\end{document}